\documentclass[10pt,twocolumn]{article} 
\usepackage{latex8}
\usepackage{times}
\usepackage{cmll}
\usepackage[T1]{fontenc}
\usepackage[latin1]{inputenc}
\usepackage{latexsym}
\usepackage{amsmath,graphics}
\usepackage{amssymb,wasysym,stmaryrd}
\usepackage{amsfonts,verbatim}
\usepackage{graphicx,lscape}
\usepackage{epsfig,color}
\usepackage{array}
\usepackage{flafter}
\usepackage[all]{xy}
\usepackage{bussproofs}
\usepackage[mathscr]{euscript}
\usepackage{mathpartir,infer}

\newif\ifcompact \compacttrue
\newif\ifultracompact \ultracompactfalse
\newif\ifsuction \suctiontrue
\newif\iftitle \titletrue

\newcommand{\TRUE}{\mathtt{true}}
\newcommand{\FALSE}{\mathtt{false}}

\newcommand{\BOOL}{\mathbb{B}}

\newcommand{\question}{\mathtt{q}}

\newcommand{\pop}{\multimap}
\newcommand{\tensor}{\otimes}
\newcommand{\synchTensor}{\odot}
\renewcommand{\with}{\&}

\newcommand{\lambdaop}[1]{\lambda^{OP}{#1}}
\newcommand{\lambdaqa}[1]{\lambda^{QA}{#1}}

\newcommand{\Acat}{\mathcal{A}}
\newcommand{\Mcat}{\mathcal{M}}
\newcommand{\Acatcont}{\mathcal{A}^{\lnot}}
\newcommand{\cat}{\mathcal{C}at}

\newcommand{\morph}[1]{\stackrel{#1}{\longrightarrow}}

\newcommand{\pos}[1]{V}
\newcommand{\mov}[1]{E}
\newcommand{\pol}[1]{\lambda}
\newcommand{\pay}[1]{\kappa}
\renewcommand{\root}[1]{\star}
\newcommand{\req}[2]{Q_{#1}({#2})}
\newcommand{\sem}[1]{[#1]}

\newcommand{\poss}[1]{V_{#1}}

\newcommand{\payy}[1]{\kappa_{#1}}
\newcommand{\roott}[1]{\star_{#1}}

\newcommand{\path}{\twoheadrightarrow}

\newcommand{\Implies}{\hspace{1em} \Longrightarrow \hspace{1em}}
\newcommand{\define}{\stackrel{\mathrm{def}}{=}}
\newcommand{\dual}[1]{#1^*}
\newcommand{\Par}{\parr}

\newcommand{\ses}{\mathcal{B}}
\newcommand{\negative}{\mathcal{B}^-}
\newcommand{\positive}{\mathcal{B}^+}
\newcommand{\pointed}{\mathcal{B}^\bullet}
\renewcommand{\neg}[1]{#1^-}

\newcommand{\C}{\mathcal{C}}
\renewcommand{\fam}[1]{Fam(#1)}
\newcommand{\family}{Fam}

\newcommand{\affine}[1]{\, !\hspace{-0.8ex}_{{\textcolor{white}{\bullet}} \hspace{-1.2ex}w} #1}
\newcommand{\affinespecial}[1]{\, !\hspace{-0.8ex}_{{\textcolor{white}{\bullet}} \hspace{-0.9ex}w} #1}
\newcommand{\affinequestionmark}[1]{\, ?\hspace{-0.9ex}_{{\textcolor{white}{\bullet}} \hspace{-1.3ex}w} #1}
\newcommand{\affinepetit}[1]{\, !\hspace{-1ex}_{{\textcolor{white}{\bullet}} \hspace{-1.4ex}w} #1}

\newcommand{\exponential}[1]{\, !\hspace{-0.8ex}_{{\textcolor{white}{\bullet}} \hspace{-0.9ex}e} #1}
\newcommand{\exponentialpetit}[1]{\, !\hspace{-1ex}_{{\textcolor{white}{\bullet}} \hspace{-1.1ex}e} #1}
\newcommand{\questionmark}[1]{\, ?\hspace{-0.9ex}_{{\textcolor{white}{\bullet}} \hspace{-1.05ex}e} #1}

\newtheorem{prop}{Proposition}
\newtheorem{lemma}{Lemma}

\newcommand{\mathin}{$$}
\newcommand{\mathout}{$$}

\newcommand{\sequation}{\begin{equation}}

\newcommand{\sitemize}{\end{itemize} \vspace{-0.2em}}

\newcommand{\sitem}{\vspace{-0.5em} \item}

\ifsuction
\newenvironment{proof}{\vspace{-1.2em}\paragraph{Proof:}}{}
\makeatletter
\renewcommand\section{\@startsection {section}{1}{\z@}%
   {10pt plus 2pt minus 2pt}{10pt plus 2pt minus 2pt} {\large\bf}}
\renewcommand\subsection{\@startsection {subsection}{2}{\z@}%
   {9pt plus 2pt minus 2pt}{9pt plus 2pt minus 2pt} {\elvbf}}
\renewcommand\paragraph{\@startsection{paragraph}{4}{\z@}%
                                    {2ex \@plus1ex \@minus.2ex}%
                                    {-1em}%
                                    {\normalfont\normalsize\bfseries}}

\makeatother

\newenvironment{porogroph}[1]{\paragraph{#1}}{}
\newenvironment{newconstruct}[1]{\porogroph{#1:}\hspace{1em}}{\newline}
\else
\newenvironment{porogroph}[1]{\paragraph{#1}}{}

\fi

\begin{document}


\title{Resource modalities in game semantics}

\author{Paul-Andr\'e Melli\`es 
\and Nicolas Tabareau
\thanks{This work has been supported by the ANR
Invariants alg\'ebriques des syst\`emes informatiques (INVAL).
Postal address: Equipe PPS, Universit\'e Paris VII, 
2 place Jussieu, Case 7014, 75251 Paris Cedex 05, FRANCE.
Email addresses: \textbf{mellies@pps.jussieu.fr}
and \textbf{tabareau@pps.jussieu.fr}}}

\maketitle
\thispagestyle{empty}


\begin{abstract}
The description of resources in game semantics
has never achieved the simplicity and precision of linear logic,
because of a misleading conception: the belief
that linear logic is more primitive than game semantics.
We advocate the contrary here: that game semantics is conceptually
more primitive than linear logic.
%
%
%
Starting from this revised point of view, we design
a categorical model of resources in game semantics,
%
and construct an arena game model where the usual notion of bracketing
is extended to multi-bracketing in order to capture
various resource policies: linear, affine and exponential.
\end{abstract}


%





\section{Introduction}
\label{sct:intro}
\paragraph{Game semantics and linear logic.}
Game semantics is the younger sibling of linear logic:
born (or reborn) at the beginning of the 1990s, in the turmoil
produced by the recent discovery of linear logic
by Jean-Yves Girard~\cite{girard:linear-logic},
it remained under its spiritual influence for a very long time.
%
%
This ascendancy of linear logic was extraordinarily
healthy and profitable in the early days.
Properly guided, game semantics developed steadily,
following the idea that every \emph{formula} of linear logic describes a \emph{game};
and that every \emph{proof} of the formula describes a
\emph{strategy} for playing on that game.
This correspondence between formulas of linear logic
and games is supported by a series of elegant and striking analogies.
%
One basic principle of linear logic is that every formula
behaves as a resource, which disappears once consumed.
In particular, a proof of the formula
$A \pop B$
is required to deduce the conclusion~$B$ by using
(or consuming) its hypothesis~$A$ exactly once.
This principle is nicely reflected in game semantics,
by the idea that playing a game is just like consuming a resource,
the game itself.
Another basic principle of linear logic is that negation
$A\mapsto \lnot A$ is involutive.
This means that every formula~$A$ is equal
(or at least isomorphic) to the formula negated twice:
\begin{equation}
\label{equation/involutive-negation}
A \hspace{1em} \cong \hspace{1em} \lnot \lnot A.
\end{equation}
Again, this principle is nicely reflected in game semantics
by the idea that negating a game~$A$ consists in permuting
the r\^oles of the two players.
Hence, negating a game twice amounts to permuting
the r\^ole of Proponent and Opponent twice,
which is just like doing nothing.
%

The connectives of linear logic are also nicely
reflected in game semantics.
For instance, the tensor product
$A\tensor B$
of two formulas~$A$ and~$B$ is suitably interpreted as
the game (or formula)~$A$ played in parallel with
the game (or formula)~$B$, where only Opponent
may switch from a component to the other one.
Similarly, the sum
$A\oplus B$
of two formulas~$A$ and~$B$ is suitably interpreted as
the game where Proponent plays the first move,
which consists in choosing between the game~$A$
and the game~$B$, before carrying on in the selected component.
Finally, the exponential modality of linear logic
$
!A
$
applied to the formula~$A$ is suitably interpreted as
the game where several copies of the game~$A$
are played in parallel, and only Opponent is allowed
(1) to switch from a copy to another one 
and (2) to open a fresh copy of the game~$A$.

%
%
%

%
What we describe here is in essence the game semantics
of linear logic defined by Andreas Blass in~\cite{blass:linear-logic}.
Simple and elegant, the model reflects the full flavour of the resource policy
of linear logic.
It is also remarkable that this game semantics is an early predecessor
to linear logic~\cite{blass:determinacy}.
%
%

\paragraph{A schism with linear logic.}
%
The destiny of game semantics has been to emancipate itself
from linear logic in the mid--1990s, in order to comply with its own designs,
inherited from denotational semantics:
%
%
%
\begin{enumerate}
\sitem the desire to interpret \emph{programs} written
in programming languages with effects (recursion, states, etc.)
and to characterise exactly their interactive behaviour
inside \emph{fully abstract} models;
\sitem the desire to understand the algebraic principles of
programming languages and effects,
using the language of category theory.
\end{enumerate}
So, a new generation of game semantics arose, propelled
by (at least) two different lines of research:
\begin{enumerate}
\sitem Samson Abramsky and Radha Jagadeesan~\cite{abramsky-jagadeesan}
noticed that the (alternating variant of the) Blass model
does not define a categorical model of linear logic.
Worse: it does not even define a category,
for lack of associativity.
Abramsky dubs this phenomenon
the \emph{Blass problem} and describes it in~\cite{abramsky03}.
\sitem Martin Hyland and Luke Ong~\cite{hyland-ong}
introduced the notion of \emph{arena game},
and characterised the interactive behaviour of programs written
in the functional language~PCF --- the simply-typed
$\lambda$-calculus with conditional test, arithmetic and recursion.
\end{enumerate}
So, the Blass problem indicates that it is difficult
to construct a (sequential) game model
of linear logic; and at about the same time,
arena games become mainstream
although they do not define a model of linear logic.
These two reasons (at least) opened a schism between
game semantics and linear logic: it suddenly became accepted
that categories of (sequential) games and strategies
would only capture \emph{fragments} of linear logic
(intuitionistic or polarised) but not the whole thing.

On the other hand, defining the resource modalities of linear logic
for game semantics requires to reunify the two schismatic subjects.
Since the disagreement started with category theory,
this reunification should occur at the categorical level.
%
%
%
We explain (in \S \ref{sct:catModel}) how to achieve this
by \emph{relaxing} the involutive negation of linear logic
into a less constrained tensorial negation.
This negation induces in turn a \emph{linear continuation}
monad, whose unit
\begin{equation}
\label{equation/usual-negation}
A \hspace{1em} \longrightarrow \hspace{1em} \lnot \lnot A
\end{equation}
refines the isomorphism~(\ref{equation/involutive-negation})
of linear logic.
Moving from an involutive to a tensorial negation means that
we replace linear logic by a more general and primitive logic --
which we call \emph{tensorial logic}.
As we will see, this shift to tensorial logic clarifies the Blass problem,
and describes the structure of arena games.
It also enables the expressions of resource modalities in game
semantics, just as it is usually done in linear logic.
However, because the presentation of modalities
may appear difficult to readers not familiar with categorical semantics,
we prefer to recall first the notion
of \emph{well-bracketing} in arena games --- and explain
how it can be reunderstood as a resource policy,
and extended to multi-bracketing.

\paragraph{Arena games.}
Recall that an \emph{arena} is defined as a forest of rooted trees,
whose nodes are called the \emph{moves} of the game.
%
%
One writes
\vspace{-.2em}
$$
m\vdash n
$$
and says that the move~$m$ \emph{enables} the move~$n$ 
when the move~$m$ is the immediate ancestor of the move~$n$
in the arena.
Every move~$m$ is assigned a polarity~$\lambdaop{(m)}\in\{-1,+1\}$.
By convention, $\lambdaop{(m)}=+1$ when the move
is Proponent, and $\lambdaop{(m)}=-1$ when it is Opponent.
Finally, one requires that the arena is alternating:
$$
m\vdash n
\Implies
\lambdaop{(m)}= -\lambdaop{(n)}
$$
and that all roots (called \emph{opening} moves) of the arena
have the same polarity.
A typical example of arena is the boolean arena~$\BOOL$:
%
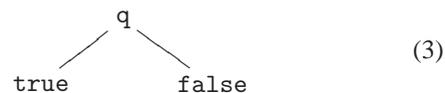
\begin{equation}
\label{equation/boolean-arena}
\vcenter{
\xymatrix @-1pc {
& \question\ar@{-}[ld]\ar@{-}[rd]
\\
\TRUE && \FALSE
}}
\end{equation}
%
where the Opponent move~$\question$ justifies
the two Proponent moves~$\TRUE$ and~$\FALSE$.
Every arena game~$A$ induces a set
of \emph{justified plays}, which are essentially sequences
of moves (we will avoid discussing \emph{pointers} here.)
%
%
Typically, the PCF type
$$(\BOOL_3\Rightarrow\BOOL_2)\Rightarrow\BOOL_1$$
defines the arena
$$
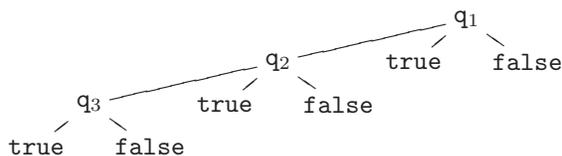
\hspace{-.3em}
\xymatrix @-2.1pc {
&&&&&&& \question_1\ar@{-}[ldd]\ar@{-}[rdd]\ar@{-}[llldd]
\\
\\
&&&& \question_2\ar@{-}[ldd]\ar@{-}[rdd]\ar@{-}[llldd]
&&\TRUE && \FALSE
\\
\\
& \question_3\ar@{-}[ldd]\ar@{-}[rdd] &&\TRUE && \FALSE
\\
\\
\TRUE && \FALSE
}$$
%
where the indices $1,2,3$ distinguish the three instances
of the boolean arena~$\BOOL$.
This arena contains the justified play
\begin{equation}
\label{equation/well-bracketed-play-flat}
\xymatrix @-2pc {
\question_1
&
\cdot
&
\question_2
&
\cdot
&
\question_3
&
\cdot
&
\TRUE_3
&
\cdot
&
\TRUE_2
&
\cdot
&
\TRUE_1
}
\end{equation}
%
%
also depicted using the convention below:
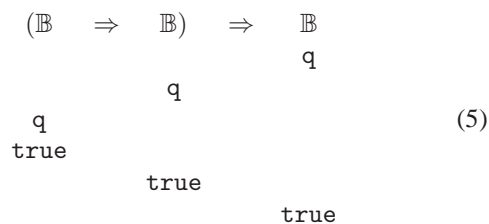
\begin{equation}\label{equation/well-bracketed-play}
\begin{array}{ccccc}
(\BOOL & \Rightarrow & \BOOL) & \Rightarrow & \BOOL
\\
&&&& \question
\\
&& \question &&
\\
\question &&&&
\\
\TRUE
\\
&& \TRUE
\\
&&&& \TRUE
\end{array}
\end{equation}
%
Note that
the play~(\ref{equation/well-bracketed-play-flat}-\ref{equation/well-bracketed-play})
belongs to the strategy implemented by the PCF program~$\lambda f.f(\TRUE)$.

\paragraph{Well-bracketing.}
Hyland and Ong demonstrate in their work~\cite{hyland-ong}
that a (finite) strategy can be implemented in PCF
if and only if it satisfies two fundamental conditions,
called \emph{innocence} and \emph{well-bracketing}.
%
%
We will focus here on the well-bracketing condition,
which is very similar to a \emph{stack discipline}.
%
%
The condition is usually expressed in the following way.
Arenas are refined by attaching
a \emph{mode}~$\lambdaqa{(m)}\in\{Q,A\}$
to every move~$m$ of the arena.
A move~$m$ is called a \emph{question} when~$\lambdaqa{(m)}=Q$,
and an \emph{answer} when~$\lambdaqa{(m)}=A$.
One then requires that no answer move~$m$
justifies another answer move~$n$:
\vspace{-.1em}
$$
m\vdash n
\Implies
\lambdaqa{(m)}=Q
\hspace{.5em}
\mbox{or}
\hspace{.5em}
\lambdaqa{(n)}=Q.
$$
%
The intuition indeed is that an answer~$n$ responds
to the question~$m$ which justifies it in the play.
Note that alternation ensures that Proponent answers
the questions raised by Opponent, and vice versa:
hence, a player never answers his own questions.
For instance, the arena game~$\BOOL$ is refined by declaring
that the Opponent move~$\question$ is a question, and
that the two Proponent moves~$\TRUE$ and~$\FALSE$ are answers.

%

%
Now, a justified play~$s$ is called \emph{well-bracketed}
when every answer~$n$ appearing in the play
responds to the ``pending'' question~$m$.
The terminology is supported by the intuition
that (1) every question ``opens'' a bracket and
(2) every answer ``closes'' a bracket,
which should match the bracket opened
by the answered question.
Typically, the play~(\ref{equation/well-bracketed-play-flat}-\ref{equation/well-bracketed-play})
is well-bracketed, because every answer responds properly
to the last unanswered question, thus leading to the well-bracketed
sequence:
\vspace{-.2em}
$$
\xymatrix @-2.1pc {
\question_1 & \cdot & \question_2 & \cdot & \question_3
& \cdot & \TRUE_3 & \cdot & \TRUE_2 & \cdot & \TRUE_1
\\
(_1 \ar@{-}[rrrrrrrrrr] 
&& && && && && _1)
\\
&& (_2  \ar@{-}[rrrrrr]
  && && && _2)
\\
&&&& (_3 \ar@{-}[rr] && _3)
}
\vspace{-.2em}
$$
%
On the other hand, the play
%
\begin{equation}\label{equation/non-well-bracketed-play}
\begin{array}{ccccc}
(\BOOL & \Rightarrow & \BOOL) & \Rightarrow & \BOOL
\\
&&&& \question
\\
&& \question &&
\\
\question &&&&
\\
&&&& \TRUE
\end{array}
\end{equation}
%
is \emph{not} well-bracketed, because the move~$\TRUE$
answers the first question of the play, whereas it should
have answered the third (and pending) question.
This may be depicted in the following way:
\begin{equation}\label{equation/non-well-bracketed-play-flat}
\vcenter{\xymatrix @-2.1pc {
\question_1 & \cdot & \question_2 & \cdot & \question_3 & \cdot & \TRUE_1
\\
(_1\ar@{-}[rrrrrr] && && && _1)
\\
&& (_2 \ar@{.}[rrrr] && && &&
\\
&& && (_3 \ar@{.}[rr] && &&
}}
\end{equation}
%
%
In fact, the
play~(\ref{equation/non-well-bracketed-play}-\ref{equation/non-well-bracketed-play-flat})
belongs to a strategy which tests whether
the function~$f:\BOOL\Rightarrow\BOOL$ is strict, that is, interrogates
its argument: this test cannot be implemented in the language PCF
-- although it can be implemented in PCF extended with the control operator
\texttt{call-cc}, see \cite{cartwright1994,laird97}.
\paragraph{Counting resources.}
We would like to understand well-bracketing as a resource discipline,
rather than simply as a stack discipline.
One key step in this direction is the observation that a
well-bracketed play may be detected simply by counting two specific
numbers on a path:
\begin{itemize}
\sitem the number~$\pay A^{+}$ of Proponent questions opened
but left unanswered,
\sitem the number~$\pay A^{-}$ of Opponent questions opened
but left unanswered.
\end{itemize}
Of course, it is not sufficient to count the two numbers~$\pay A^+$
and~$\pay A^-$ of a play~$s$ to detect whether the play is well-bracketed.
Typically, the well-bracketed play~$(a)$ and the non well-bracketed play~$(b)$
introduced in~(\ref{equation/non-well-bracketed-play}-\ref{equation/non-well-bracketed-play-flat}) induce the same numbers $\pay A^+$ and $\pay A^-$:
$$
\begin{array}{cccccc}
(a) & &
\question_1 \cdot \question_2 \cdot \question_3 \cdot \TRUE_3
&
\longmapsto
&
\pay A^+=1,
&
\pay A^-=1 
\\
(b) &&
\question_1 \cdot \question_2 \cdot \question_3 \cdot \TRUE_1
&
\longmapsto
&
\pay A^+=1, 
&
\pay A^-=1 
\end{array}
$$
In order to detect well-bracketing, one needs to apply
the count to the subpaths $(c)$ and $(d)$ of these plays.
This reveals a key difference:
$$
\begin{array}{cccccc}
(c) &&
\question_3 \cdot \TRUE_3
&
\longmapsto
&
\pay A^+=0, 
&
\pay A^-=0
\\
(d)
&&
\question_3 \cdot \TRUE_1
&
\longmapsto
&
\pay A^+=0,
&
\pay A^-=1 
\end{array}
$$
The elementary but key characterisation follows:
\begin{prop}
\label{proposition/well-bracketing}
A play~$s$ is well-bracketed if and only if every subpath~$m\cdot t\cdot n$
of the play~$s$ satisfies
$$
\pay A^+(m\cdot t\cdot n) = 0 \Implies
\pay A^-(m\cdot t\cdot n) = 0
$$
when $m$ is Opponent and $n$ is Proponent;
and dually
$$
\pay A^-(m\cdot t\cdot n) = 0 \Implies
\pay A^+(m\cdot t\cdot n) = 0
$$
when $m$ is Proponent and $n$ is Opponent.
\end{prop}
Let us explain this briefly.
Suppose that~$m\cdot t\cdot n$ is a subpath of a well-bracketed play~$s$,
where~$m$ is Opponent and $n$ is Proponent.
The first condition says that if there is an Opponent question
unanswered in~$m\cdot t$, then either Player answers it
-- in which case $\pay A^-(m\cdot t\cdot n) = 0$ --
or there is a Player question unanswered in~$m\cdot t\cdot n$
-- in which case $\pay A^+(m\cdot t\cdot n) \neq 0$.
The other condition is dual.

\paragraph{A resource policy.}
Reformulated in this way, the well-bracketing looks very much like a resource policy.
%
%
The basic intuition is that every question~$m$ emits
a  \emph{query} for a \emph{linear session}.
This query is noted by a opening bracket~$(_i$ and
counted by~$\pay{}^{±}$ where $±$ is the polarity
of the move~$m$.
The query is then complied with by a \emph{response}
emitted by an answer move~$n$, and noted
by a closing bracket~$_i )$.
In our example, the move~$\question_3$ emits a query~$(_3$
which is later complied with
in the play~(\ref{equation/well-bracketed-play-flat}-\ref{equation/well-bracketed-play})
by the response~$_3 )$ emitted by the move~$\TRUE$ whereas
it remains unanswered in the
play~(\ref{equation/non-well-bracketed-play}-\ref{equation/non-well-bracketed-play-flat}).
Hence, a play like~(\ref{equation/non-well-bracketed-play}-\ref{equation/non-well-bracketed-play-flat}) is not well-bracketed because it breaks the linearity policy
implemented by the queries.
Our game model will relate this linearity policy to the fact
that the boolean formula is defined as
\begin{equation}
\label{equation/linear-boolean}
\BOOL \hspace{1em} = \hspace{1em} 
\stackrel{O}{\lnot}\ \stackrel{P}{\lnot} \ (1\oplus 1)
\end{equation}
in tensorial logic.
Here, the tags~$O$ and~$P$ are mnemonics
to indicate that the external negation~$\lnot_O$
is interpreted as an Opponent move, whereas the internal negation~$\lnot_P$
is interpreted as a Proponent move.
The story told by~(\ref{equation/linear-boolean}) goes like this:
Opponent plays the external negation, followed by Proponent,
who plays the internal negation and \emph{at the same time}
resolves the choice $1\oplus 1$ between $\TRUE$ and $\FALSE$.
This refines the picture conveyed by the boolean
arena~(\ref{equation/boolean-arena}) by decomposing the Player
moves~$\TRUE$ and~$\FALSE$ in two compound stages: negation and choice
-- where negation thus encapsulates the two moves~$\TRUE$ and~$\FALSE$.
This enables to relax the well-bracketing policy by interpreting
the boolean formula as
\begin{equation}
\label{equation/affine-boolean}
\BOOL \hspace{1em} = \hspace{1em} 
\stackrel{O}{\lnot}\ \ \  \affine{} \ \  \stackrel{P}{\lnot} \ (1\oplus 1)
\end{equation}
where the affine modality~$\affine{}$ of tensorial logic
is inserted between the two negations.
The intuitionistic hierarchy on the boolean
formula~(\ref{equation/linear-boolean}) coincides with the
well-bracketed arena game model of PCF described by Hyland and
Ong in~\cite{hyland-ong} whereas the intuitionistic hierarchy on the
boolean formula~(\ref{equation/affine-boolean}) -- where the affine
modality $\affine{}$ is replaced by the exponential modality
$\exponential{}$ -- coincides with the non-well-bracketed arena game
model of PCF with control described by Jim Laird in~\cite{laird97}
and Olivier Laurent in~\cite{laurent:games}.

%
%

\paragraph{Multi-bracketing.}
This analysis leads us to the notion of \emph{multi-bracketing} in arena games.
In linear logic, every proof of the formula
$$(\BOOL \tensor \BOOL) \pop \BOOL$$
asks the value of its two boolean arguments,
and we would like to understand this as a kind
of well-bracketing condition.
%
So, the play
%
\begin{equation}
\label{equation/generalized-well-bracketed}
\begin{array}{ccccc}
(\BOOL & \tensor & \BOOL) & \pop & \BOOL
\\
&&&& \question
\\
&& \question &&
\\
&& \TRUE &&
\\
\question &&&&
\\
\TRUE &&&&
\\
&&&& \TRUE
\end{array}
\end{equation}
%
would be ``well-bracketed'' in the new setting,
whereas the play
%
\begin{equation}
\label{equation/non-generalized-well-bracketed}
\begin{array}{ccccc}
(\BOOL & \tensor & \BOOL) & \pop & \BOOL
\\
&&&& \question
\\
\question &&&&
\\
\TRUE &&&&
\\
&&&& \TRUE
\end{array}
\end{equation}
%
%
would not be ``well-bracketed'',
because it does not explore the second argument of the function.
This extended well-bracketing is captured by the idea
that the first question emits \emph{three}
queries~$(_1$ and $(_a$ and $(_b$ at the same time.
Then, the play~(\ref{equation/generalized-well-bracketed})
appears to be ``well-bracketed'' if one depicts the situation
in the following way:
$$
\xymatrix @-2.1pc {
\question_1 & \cdot & \question_2 & \cdot & \TRUE_2
& \cdot & \question_3 & \cdot & \TRUE_3 & \cdot & \TRUE_1
\\
(_{1} \ar@{-}[rrrrrrrrrr] &&&&&&&&&& _{1})
\\
(_{a} \ar@{-}[rr] && _{a}) (_2 \ar@{-}[rr] && _2) && && &&
\\
(_{b} \ar@{-}[rrrrrr] && && && _{b}) (_3 \ar@{-}[rr] && _3) &&
}
$$
whereas the play~(\ref{equation/non-generalized-well-bracketed})
is not ``well-bracketed'' because the query~$(_a$
is never complied with, as can be guessed from the picture below:
$$
\xymatrix @-2.1pc {
\question_1 & \cdot & \question_3 & \cdot & \TRUE_3
& \cdot & \TRUE_1
\\
(_{1} \ar@{-}[rrrrrr] &&&&&& _{1})
\\
(_{a} \ar@{.}[rrrrrr] && && && 
\\
(_{b} \ar@{-}[rr] && _{b}) (_3 \ar@{-}[rr] && _3) && && &&
}
$$
We explain in $\S \ref{sct:conway games}$ and $\S \ref{sct:game}$ how we apply
the well-bracketing criterion devised in
Proposition~\ref{proposition/well-bracketing} in order to generalise
well-bracketing to a multi-bracketed framework.

\paragraph{Plan of the paper.}
We describe ($\S \ref{sct:catModel}$) a categorical semantics of resources
in game semantics, and explain in what sense the resulting topography
refines both linear logic and polarized logic.
After that, we construct ($\S \ref{sct:conway games}$) a compact-closed (that is, self-dual)
category of multi-bracketed Conway games
and well-bracketed strategies,
where the resource policy is enforced by multi-bracketing.
From this, we derive ($\S \ref{sct:game}$) a model of our categorical
semantics of resources, using a family construction, and conclude ($\S
\ref{sct:conclusion}$).

\paragraph{Acknowledgements.}
We would like to thank Martin Hyland together with Masahito Hasegawa,
Olivier Laurent, Laurent Regnier and Peter Selinger for stimulating
discussions at various stages of this work.

\section{Categorical models of resources}\label{sct:catModel}
We introduce now the notion of \emph{tensorial negation}
on a symmetric monoidal category; and then explain
how such a category with negation may be equipped
with additives and various resource modalities.
%
%
%
The first author describes in~\cite{mellies:functorial-boxes}
how to extract a \emph{syntax} of proofs from a categorical semantics,
using string diagrams and functorial boxes.
The recipe may be applied here to extract the syntax of a logic,
called \emph{tensorial logic.}
However, we provide in Appendix  a sequent calculus
for tensorial logic, in order to compare it to linear logic~\cite{girard:linear-logic}
or polarized linear logic~\cite{laurent:phd}.

\newcommand{\llunit}{1}

\paragraph{Tensorial negation.}
A \emph{tensorial negation} on a symmetric monoidal
category~$(\Acat,\tensor,\llunit)$ is defined as a functor
\begin{center}
\begin{tabular}{ccccc}
$\lnot$ & : & $\Acat$ & $\longrightarrow$ & $\Acat^{op}$
\end{tabular}
\end{center}
together with a family of bijections
$$
\begin{array}{ccccc}
\varphi_{A,B,C}
&
:
&
\Acat(A\tensor B,\lnot{C})
&
\cong
&
\Acat(A,\lnot{(B\tensor C)})
\end{array}
$$
natural in~$A,B$ and~$C$.
Given a negation, it is customary to define
the formula \emph{false} as the object
$$
\bot \hspace{.2em} \define \hspace{.2em} \lnot \ \llunit
$$
obtained by ``negating'' the unit object~$\llunit$ of the monoidal category.
Note that we use the notation~$1$ (instead of $I$ or $e$)
in order to remain consistent with the notations of linear logic.
Note also that the bijection~$\varphi_{A,B,\llunit}$ provides then
the category~$\Acat$ with a one-to-one correspondence
$$
\begin{array}{ccccc}
\varphi_{A,B,\llunit}
&
:
&
\Acat(A\tensor B,\bot)
&
\cong
&
\Acat(A,\lnot B)
\end{array}
$$
for all objects~$A$ and $B$.
For that reason, the definition of a negation~$\lnot$ is often replaced
by the --- somewhat too informal --- statement that
``the object~$\bot$ is exponentiable'' in the symmetric
monoidal category~$\Acat$, with negation $\lnot A$ noted~$\bot^{A}$.

\paragraph{Self-adjunction.}
In his PhD thesis, Hayo Thielecke~\cite{thielecke:thesis}
observes for the first time a fundamental ``self-adjunction''
phenomenon, related to negation.
This observation plays then a key r\^ole in an unpublished work
by Peter Selinger and the first author~\cite{mellies-selinger}
on polar categories, a categorical semantics
of polarized linear logic, continuations and games.
The same idea reappears recently in a nice, comprehensive study
on polarized categories (=distributors) by Robin Cockett
and Robert Seely~\cite{cockett-seely:games}.
%
%
In our situation, the ``self-adjunction'' phenomenon amounts
to the fact that every tensorial negation is left adjoint
to the opposite functor
%
\begin{equation}
\label{equation/negation-op}
\begin{array}{ccccc}
\lnot & : & \Acat^{op} & \longrightarrow & \Acat
\end{array}
\end{equation}
because of the natural bijection
$$
\Acat^{op}(\lnot A, B)
\hspace{1em}
\cong
\hspace{1em}
\Acat(A, \lnot B).
$$

\paragraph{Continuation monad.}
Every tensorial negation~$\lnot$ induces an adjunction,
and thus a monad
$$
\lnot\lnot \hspace{1em}: \hspace{1em}
\Acat \longrightarrow \Acat
$$
This monad is called the \emph{continuation monad}
of the negation.
One fundamental fact observed by Eugenio Moggi~\cite{moggi}
is that the continuation monad is \emph{strong}
but not commutative in general.
By strong monad, we mean that the monad~$\lnot\lnot$
is equipped with a family of morphisms:
$$
t_{A,B}
\hspace{1em}
:
\hspace{1em}
A\tensor \lnot\lnot B \hspace{1em}
\longrightarrow 
\hspace{1em}
\lnot\lnot \ (A\tensor B)
$$
natural in~$A$ and~$B$, and satisfying a series
of coherence properties.
By commutative monad, we mean a strong monad
making the two canonical morphisms
%
\begin{equation}
\label{equation/two-arrows}
\lnot \lnot A \tensor \lnot\lnot B 
\hspace{1em}
\rightrightarrows
\hspace{1em}
\lnot\lnot \ (A\tensor B)
\end{equation}
coincide.
A tensorial negation~$\lnot$ is called \emph{commutative}
when the continuation monad induced in~$\Acat$ is commutative
--- or equivalently, a monoidal monad in the lax sense.

\paragraph{Linear implication.}
A symmetric monoidal category~$\Acat$
with a tensorial negation~$\lnot$ is not very far
from being monoidal \emph{closed.}
It is possible indeed to define a \emph{linear implication}~$\pop$
when its target~$\lnot B$ is a negated object:
\vspace{-.4em}
$$
\vspace{-.4em}
A\ \pop \ \lnot B \quad \define \quad
\lnot \ (A\tensor B).
$$
In this way, the functor~(\ref{equation/negation-op}) defines
what we call an \emph{exponential ideal} in the category~$\Acat$.
When the functor is faithful on objects and morphisms,
we may identify this exponential ideal with the subcategory
of \emph{negated objects} in the category~$\Acat$.
The exponential ideal discussed in Guy McCusker's PhD
thesis~\cite{mccusker:phd} arises precisely in this way.
%
This enables in particular to define the linear
and intuitionistic hierarchies
on the arena games~(\ref{equation/linear-boolean})
and~(\ref{equation/affine-boolean}).

\paragraph{Continuation category.}
Every symmetric monoidal category~$\Acat$
equipped  with a negation~$\lnot$ induces
a \emph{category of continuations}~$\Acatcont$
with the same objects as~$\Acat$, and morphisms defined as
\vspace{-.4em}
$$
\vspace{-.4em}
 \Acatcont (A,B) \quad \define \quad
\Acat (\lnot A,\lnot B).
$$
Note that the category~$\Acatcont$ is the kleisli category
associated to the comonad in~$\Acat^{op}$ induced by the adjunction;
and that it is at the same time the opposite of the kleisli category
associated to the continuation monad in $\Acat$.
Because the continuation monad is strong,
the category~$\Acatcont$ is \emph{premonoidal} in the sense
of John Power and Edmund Robinson~\cite{power-robinson}.
It should be noted that string diagrams in premonoidal categories are inherently
related to control flow charts in software engineering, as noticed
by Alan Jeffrey~\cite{jeffrey:premonoidal}.

\paragraph{Semantics of resources.}
A \emph{resource modality} on a symmetric monoidal category~$(\Acat,\tensor,e)$
is defined as an adjunction:
\vspace{-.6em}
\begin{equation}
\label{equation/adjunction}
\xymatrix{\Mcat\ar@/^1pc/[rr]^{U}
\ar@{}[rr]|-{\bot}&& \Acat\ar@/^1pc/[ll]^{F}}
\vspace{-1em}
\end{equation}
where
\begin{itemize}
\sitem $(\Mcat,\bullet,u)$ is a symmetric monoidal category,
\sitem $U$ is a symmetric monoidal functor.
\end{itemize}
Recall that a \emph{symmetric monoidal} functor~$U$
is a functor which transports the symmetric monoidal structure
of~$(\Mcat,\bullet,u)$ to the symmetric monoidal structure of~$(\Acat,\tensor,e)$,
up to isomorphisms satisfying suitable coherence properties.
Another more conceptual definition of a resource modality
is possible: it is an adjunction defined in the 2-category of symmetric monoidal categories,
\emph{lax} symmetric monoidal functors, and monoidal transformations.
Now, the resource modality is called
\begin{itemize}
\sitem \emph{affine} when the unit~$u$ is the terminal object
of the category~$\Mcat$,
\sitem \emph{exponential} when the tensor product~$\bullet$ is 
a cartesian product,
and the unit~$u$ is the terminal object of the category~$\Mcat$.
\end{itemize}
This definition of resource modality is inspired by the categorical
semantics of linear logic, and more specifically by Nick Benton's notion
of Linear-Non-Linear model \cite{benton:linear_non_linear_logic} ---
which may be reformulated now
as a symmetric monoidal \emph{closed} category~$\Acat$
equipped with an exponential modality
in our sense.
Very often, we will identify the resource modality
and the induced comonad~$! \ = \ U\circ F$ on the category~$\Acat$.
%

\paragraph{Tensorial logic.}
In our philosophy, tensorial logic is entirely described
by its categorical semantics --- which is defined
in the following way.
First, every symmetric monoidal category~$\Acat$
equipped with a tensorial negation~$\lnot$
defines a model of \emph{multiplicative} tensorial logic.
%
%
Such a category defines a model of \emph{multiplicative additive}
tensorial logic when the category~$\Acat$ has finite coproducts (noted~$\oplus$)
which \emph{distribute} over the tensor product: this means that the canonical morphisms
\vspace{-.8em}
$$
\vspace{-.5em}
(A\tensor B) \oplus (A \tensor C)
\hspace{1em}
\morph{}
\hspace{1em}
A \tensor (B \oplus C)
$$
\vspace{-1.5em}
$$
\vspace{-.5em}
0
\hspace{1em}
\morph{}
\hspace{1em}
A \tensor 0
$$
are isomorphisms.
Then, a model of (full) tensorial logic is defined as
a model of multiplicative additive tensorial logic,
equipped with an affine resource modality
(with comonad noted~$\affine$)
as well as an exponential resource modality
(with comonad noted~$\exponential $).

The diagrammatic syntax of tensorial logic will be
readily extracted from its categorical definition, using
the recipe explained in~\cite{mellies:functorial-boxes}.
However, the reader will find a sequent calculus of tensorial logic
in Appendix, written in the more familiar fashion of proof theory.
Seen from that point of view, the modality-free fragment of tensorial logic
describes a \emph{linear} variant of Girard's LC~\cite{girard:lc}
thus akin to ludics~\cite{girard:ludics} and more precisely
to what Laurent calls MALLP in his PhD thesis~\cite{laurent:phd}.
This convergence simply expresses the fact that these systems
are all based on tensors, sums and linear continuations.

%

%

\paragraph{Arena games and classical logic.}
Starting from Thielecke's work, Selinger~\cite{selinger:control}
designs the notion of \emph{control category} in order
to axiomatize the categorical semantics of classical logic.
Then, prompted by a completeness result established
by Martin Hofmann and Thomas Streicher
in~\cite{hofmann-streicher02}, he proves
a beautiful structure theorem,
stating that every control category~$\C$ is the continuation
category~$\Acatcont$ of a \emph{response category}~$\Acat$.
Now, a response category~$\Acat$
--- where the monic requirement on
the units~(\ref{equation/usual-negation}) is relaxed ---
is exactly the same thing as a model of multiplicative additive
tensorial logic, where the tensor~$\tensor$ is \emph{cartesian}
and the tensor unit~$1$ is \emph{terminal.}
A purely proof-theoretic analysis of classical logic leads exactly
to the same conclusion.
Starting from Girard's work on polarities in LC~\cite{girard:lc} and
ludics \cite{girard:ludics}, Laurent developed a comprehensive analysis
of polarities in logic, incorporating classical logic, control categories
and (non-well-bracketed) arena games~\cite{laurent:phd, laurent:games}.
Now, it appears that Laurent's polarized logic LLP coincides
with multiplicative additive tensorial logic
--- where the monoidal structure is \emph{cartesian}.
%
This is manifest in the monolateral formulation
of tensorial logic, see Appendix.
We sum up below the difference between tensorial logic
and classical logic in a very schematic table:
\begin{center}
\begin{tabular}{|c|c|}
\hline
Tensorial logic &
\begin{tabular}{cc}
$\tensor$ & is monoidal
\\
$\lnot$ & is tensorial
\end{tabular}
\\
\hline
Classical logic &
\begin{tabular}{cc}
$\tensor$ & is cartesian
\\
$\lnot$ & is tensorial
\end{tabular}
\\
\hline
\end{tabular}
\end{center}
Note that every resource modality~(\ref{equation/adjunction})
on a category~$\Acat$ equipped with a tensorial negation~$\lnot$ induces
a tensorial negation~$F^{op}\circ\lnot\circ U$ on the category~$\Mcat$.
This provides a model of polarized linear logic, and thus of classical logic,
whenever~$\Mcat$ is cartesian.
This phenomenon underlies the construction of a control category
in~\cite{LaurentRegnier03}, see also~\cite{hamano2007}
for another construction.
%

\paragraph{Linear logic.}
The continuation
monad~$A \hspace {.1em} \mapsto \hspace{.5em} \stackrel{O}{\lnot}\stackrel{P}{\lnot} A$
of game semantics lifts an Opponent-starting game~$A$ with an Opponent move~$\lnot_O$
followed by a Player move~$\lnot_P$.
Now, it appears that the Blass problem mentioned in $\S \ref{sct:intro}$ arises
precisely from the fact that the monad is strong,
but not commutative~\cite{mellies-selinger,mellies:ag3}.
%
%
%
Indeed, one obtains a game model of (full) propositional linear logic
by~\emph{identifying} the two canonical strategies~(\ref{equation/two-arrows})
--- this leading to a fully complete model of linear logic expressed
in the language of asynchronous games~\cite{mellies:ag4}.
This construction in game semantics has a nice categorical counterpart.
We already mentioned that the continuation category~$\Acatcont$
inherits a premonoidal structure from the symmetric monoidal
structure of~$\Acat$.
Now, Hasegawa Masahito shows (private communication)
that 
the continuation category~$\Acatcont$
equipped with this premonoidal structure
is $\ast$-autonomous if and only if
the continuation monad
is commutative.
The specialist will recognize here a categorification
of Girard's phase space semantics~\cite{girard:linear-logic}.
Anyway, this shows that linear logic is essentially tensorial logic 
in which the tensorial negation is commutative.
\begin{center}
\begin{tabular}{|c|c|}
\hline
Linear logic &
\begin{tabular}{cc}
$\tensor$ & is monoidal
\\
$\lnot$ & is commutative
\end{tabular}
\\
\hline
\end{tabular}
\end{center}
In that situation, every resource modality on the category~$\Acat$
induces a resource modality on the $\ast$-autonomous
category~$\Acatcont$, and thus a model of full linear logic.
%


\section{Multi-bracketed Conway games} \label{sct:conway games}
We define here and in $\S \ref{sct:game}$
a game semantics with resource modalities and fixpoints,
in order to interpret recursion in programming languages.
We achieve this by constructing first a compact-closed category~$\ses$
of \emph{multi-bracketed} Conway games, inspired from Andr\'e Joyal's
pioneering work~\cite{joyal}.
The compact-closed structure of~$\ses$ induces a trace
operator~\cite{joyal-street-verity} which, in turn, provides enough
fixpoints in the category constructed in~$\S \ref{sct:game}$ in order
to interpret the language~PCF enriched with resource modalities.

%
%
\paragraph{{\bf Conway games.}}
A Conway game is an oriented rooted graph $(\pos A,\mov A,\pol A)$
consisting of a set $\pos A$ of vertices called the \emph{positions}
of the game, a set $\mov A \subset \pos A \times \pos A$ of edges
called the \emph{moves} of the game, a function $\pol A: \mov A
\rightarrow \{-1,+1\}$ indicating whether a move belongs
to Opponent ($-1$) or Proponent ($+1$).
We note $\root A$ the root of the underlying graph.
\paragraph{{\bf Path and play.}} 
A \emph{play} is a path starting from the root $\star_A$ of the
multi-bracketed game: 
\begin{equation}
\label{equation/play}
\star_A \xrightarrow{m_1} x_1 \xrightarrow{m_2} \ldots
\xrightarrow{m_{k-1}} x_{k-1} \xrightarrow{m_k} x_k 
\end{equation}
Two paths are parallel when they have the same initial and
final positions.
A play~(\ref{equation/play}) is \emph{alternating} when:
\mathin
\forall i \in \{1 , \ldots , k-1\}, \quad \quad \lambda_A(m_{i+1}) = - \lambda_A(m_i).
\mathout 
\paragraph{{\bf Strategy.}}
A strategy $\sigma$ of a Conway game is defined as
a set of alternating plays of even length such that:
\begin{itemize}
\sitem
  $\sigma$~contains the empty play,
\sitem
  every nonempty play starts with an Opponent move,
\sitem
  $\sigma$ is closed by even-length prefix:
  for every play~$s$, and for all moves~$m,n$,
  $$
  s \cdot m \cdot n \in \sigma \Implies s \in \sigma,
  $$
\item
  $\sigma$ is deterministic: for every play~$s$,
  and for all moves~$m,n,n'$,
  \mathin
  s \cdot m \cdot n \in \sigma \mbox{ and } s \cdot m \cdot n' \in
  \sigma \ \Implies \ n=n'.
  \mathout
\sitemize
We write $\sigma: A$ when $\sigma$ is a strategy of $A$.
Note that a play in a Conway game is generally non-alternating,
but that alternation is required on the plays of a strategy.
\paragraph{Multi-bracketed games.}
A multi-bracketed game is a Conway game equipped with
\begin{itemize}
\sitem
  a finite set~$\req{A}{x}$ of \emph{queries} for each position~${x \in \pos A}$
  of the game,
\sitem
  a function $\pol{A}(x):\req{A}{x}\morph{}\{-1,+1\}$ which assigns
  to every query in $\req{A}{x}$ a polarity
  which indicates whether the query is made by Opponent ($-1$)
  or Proponent ($+1$),
  \sitem
  for each move $x \xrightarrow{m} y$, a \emph{residual} relation
    $$\sem m \subset \req{A}{x} \times \req{A}{y}
$$
  satisfying:
  \mathin
  \begin{array}{rcl}
    r \sem m r_1 \hspace{1em} \mbox{and} \hspace{1em} 
    r \sem m r_2 & \implies & r_1 = r_2 \\  
    r_1 \sem m r \hspace{1em} \mbox{and} \hspace{1em}
    r_2 \sem m r & \implies & r_1 = r_2  
  \end{array}
  \mathout 
\sitemize
The definition of residuals is then extended
to paths~$s:x\path y$ in the usual way:
by composition of relations.
We then define
  $$
    r[s]\define\{r' \ | \ r[s]r'\}
\quad \mbox{and} \quad
[s]r\define\{r' \ | \ r'[s]r\}.$$
We say that a path~$s:x\path y$:
\begin{itemize}
\item \emph{complies with a query}~$r\in\req{A}{x}$
when~$r$ has no residual after $s$ --- that is, $r[s]=\emptyset$,
\vspace{-.5em}
\item \emph{initiates a query}~$r\in\req{A}{y}$ when~$r$
has no ancestor before~$s$ --- that is, $[s]r=\emptyset$.  
\end{itemize}
We require that a move~$m$ only initiates queries of its own polarity,
and only complies with queries of the opposite polarity.
In order to formalise that a residual of a query is intuitively the query
itself, we also require that two parallel paths $s$ and $t$
induce the same residual relation: $\sem{s}{} = \sem{t}{}$.
Finally, we require that there are no queries at the root:
$\req{A}{\roott{}}=\emptyset$.

\paragraph{Resource function.}
Extending Conway games with queries enables the definition of a
resource function
\vspace{-.5em}
$$
\vspace{-.5em}
\pay A=(\pay A^+,\pay A^-)$$
which counts, for
every path $s:x\path y$, the number~$\pay A^+(s)$ (respectively~$\pay
A^-(s)$) of Proponent (respectively Opponent) queries
in~$r\in\req{A}{y}$ initiated by the path~$s$ --- that is, such that
$[s]r=\emptyset$.
The definition of multi-bracketed games induces
three cardinal properties of~$\kappa^{\pm}$, which
will replace the very definition of~$\kappa$,
and will play the r\^ole of axioms in all our proofs --
in particular, in the proof that the composite of
two well-bracketed strategies is also well-bracketed.

\smallskip \noindent 
{\bf \emph{Property 1: accuracy.}}  
For all paths~$s:x\path y$ and Proponent move~$m:y\rightarrow z,$
\vspace{-.5em}
$$
\vspace{-.5em}
\pay A^-(m) = 0 \hspace{1em} \mbox{and} \hspace{1em}
\pay a ^+(s \cdot m) = \pay A^+(s) + \pay A^+(m),
$$
as well as the dual equalities for Opponent moves.


\smallskip \noindent
{\bf \emph{Property 2: suffix domination.}}
For all paths $s:x\path y$ and $t:y \path z$,
\vspace{-.5em}
$$
\vspace{-.5em}
\pay A(t) \leq \pay A(s\cdot t).
$$

\smallskip \noindent
{\bf \emph{Property 3: sub-additivity.}}
For all paths $s:x\path y$ and $t:y \path z$,
\vspace{-.5em}
$$
\vspace{-.5em}
\pay A(s \cdot t) \leq
\pay A(s) + \pay A(t).
$$
Accuracy holds because Player does not initiate Opponent queries,
and does not comply with Player queries.
Suffix domination says that a query cannot already have been
complied with.
Sub-additivity expresses that composing two paths does
not increase the number of queries.

\paragraph{Well-bracketed plays and strategies.}

%
Once the resource function~$\pay{}$ is defined on paths,
it becomes possible to define a \emph{well-bracketed play}
as a play which satisfies the two conditions stated
in Proposition~\ref{proposition/well-bracketing} of $\S \ref{sct:intro}$.
So, the property becomes a definition here.
A strategy~$\sigma$ is then declared \emph{well-bracketed}
when, for every play~$s\cdot m \cdot t \cdot n$ of the strategy~$\sigma$
where $m$ is an Opponent move and $n$ is (necessarily)
a Proponent move:
\vspace{-.4em}
$$
\vspace{-.4em}
\payy A^+(m \cdot t \cdot n) = 0 \Implies 
\payy A^-(m \cdot t \cdot n) = 0.
$$
Every well-bracketed strategy~$\sigma$ then preserves well-bracketing
in the following sense:
%
%
\begin{lemma}
 Suppose $s \cdot m \cdot n  \in \sigma$ and that $s\cdot m$ is
 well-bracketed. Then, $s \cdot m \cdot n$ is well-bracketed.
\end{lemma}
Hence, when Opponent and Proponent play according to well-bracketed
strategies, the resulting play is well-bracketed.

%
%
%
%

%
\paragraph{Dual.}
Every multi-bracketed game~$A$ induces
a \emph{dual} game $\dual A$ obtained by reversing
the polarity of moves and queries. Thus,
$(\payy{\dual A}^+,\payy{\dual A}^-) = (\payy A^-,\payy A^+)$.

\paragraph{{\bf Tensor product.}}

The tensor product~$A\tensor B$ of two multi-bracketed games $A$ and $B$
is defined as:
\begin{itemize}
\sitem[-]
  its positions are the pairs $(x,y)$ noted $x \tensor y$, ie.  $\poss{A
  \tensor B} = \poss A \times \poss B$ with $\roott{A \tensor B} =
  (\roott A,\roott B)$. 
\sitem[-]
    its moves are of two kinds:
 $$
  x \tensor y \rightarrow 
  \left\{
  \begin{array}{l}
    z \tensor y \mbox{ if } x\rightarrow z \mbox{ in the game } A, \\
    x \tensor z \mbox{ if } y\rightarrow z \mbox{ in the game } B. \\
  \end{array}
  \right.
 $$
 \vspace{-.7em}
\sitem[-] its queries at position~$x\tensor y$ are the queries
at position~$x$ in the game~$A$ and the
queries at position~$y$ in the game~$B$:
$
 \req{A\tensor B}{x \tensor y} = \req{A}{x} \uplus \req{B}{y}.$
 \end{itemize}
The polarities of moves and queries in the game~$A\tensor B$
are inherited from the games~$A$ and~$B$, and the residual relation
of a move~$m$ in the game~$A\tensor B$ is defined just
in the expected (pointwise) way.
%
%
%
The unique multi-bracketed game~$1$ with $\{ \star \}$
as underlying Conway game is the neutral element
of the tensor product.
As usual in game semantics, every play~$s$
in the game~$A \tensor B$ may be seen
as the interleaving of a play $s_{|A}$ in the game~$A$
and a play~$s_{|B}$ in the game~$B$.
More interestingly, the resource function~$\kappa$ is ``tensorial''
in the following sense:
\vspace{-.2em}
$$
\vspace{-.4em}
\payy{A \tensor B} (s) = \payy A (s_{|A}) + \payy B (s_{|B}).
$$
\vspace{-.8em}

\vspace{-.8em}
\paragraph{{\bf Composition.}}
We proceed as in~\cite{mccusker:phd,harmer:phd}, and say that
$u$ is an interaction on three games~$A,B,C$, this noted $u \in int_{ABC}$, 
when the projection of~$u$ on each game $\dual A \tensor B$, $\dual B
\tensor C$ and $\dual A \tensor C$ is a play.  Given two strategies
$\sigma: \dual A \tensor B,\tau : \dual B \tensor C$, we define the
composition of these strategies as follows:
\mathin
\sigma;\tau = \{ u_{|\dual A \tensor C} \ | \ u \in int_{ABC} ,
u_{|\dual A \tensor B} \in \sigma , u_{|\dual B \tensor C} \in \tau\} 
\mathout
As usually, the composition of two strategies is a strategy.
More interestingly, we show that our notion of well-bracketing is preserved
by composition:

\begin{prop}
  The strategy $\sigma;\tau:\dual A \tensor C$ is well-bracketed
  when the two strategies $\sigma:\dual A \tensor B$ and
  $\tau:\dual B \tensor C$ are well-bracketed.
  \end{prop}

\begin{proof}
  \hspace{-0.8em}
  The proof is entirely based on the three cardinal properties of~$\kappa$
  mentioned earlier. The proof appears in the Master's thesis of the second
  author~\cite{tabareau05}.
\end{proof}

\paragraph{{\bf The category $\ses$ of multi-bracketed games.}}
The category $\ses$ has
multi-bracketed games as objects,
and well-bracketed strategies~$\sigma$ of~$\dual A \tensor B$
as morphisms $\sigma:A \rightarrow B$.
%
The identity strategy is the usual copycat strategy,
defined by Andr\'e Joyal in Conway game~\cite{joyal}.
The resulting category~$\ses$ is compact-closed in the sense
of~\cite{kelly:compactClosed} and thus admits 
a canonical trace operator, unique up to equivalence,
see~\cite{joyal-street-verity} for details.
%
%
%
%

\paragraph{Negative and positive games.}
A multi-bracketed game~$A$ is called \emph{negative} when all the moves
starting from the root~$\roott{A}$ are Opponent moves; and \emph{positive}
when its dual game~$\dual{A}$ is negative.
%
The full subcategory of negative (resp. positive) multi-bracketed games
is noted~$\negative$ (resp. $\positive$).
For a multi-bracketed game $A$, we write $\neg{A}$ for
the negative game obtained by removing all the Player moves
from the root.

%
%
%

\paragraph{The exponential modality.}
%
%
Every multi-bracketed game~$A$ induces an \emph{exponential} game~$! A$
as follows:
%
\begin{itemize}
\sitem[-]
  its positions are the words $w=x_1\cdots x_k$
  whose letters are positions~$x_i$ of the game~$A$ different from the
  root $\roott{A}$;
  the intuition is that the letter~$x_i$ describes
  the current position of the $i^{th}$ copy of the game,
\sitem[-]
  its root~$\roott{! A}$ is the empty word,
\sitem[-]
  its moves $w\rightarrow w'$ are either moves played in one copy:
\vskip -2em
  $$
  w_1 \ x \ w_2 \hspace{.1em} \to \hspace{.1em} \ w_1 \ y \ w_2
$$
  where $x\to y$ is a move in the game~$A$;
  or moves where Opponent opens a new copy:
  \vskip -1.3em
  $$w \ \to \ w \ x$$
  \vskip -0.7em
  where $\roott{A}\to x$ is an Opponent move in~$A$.
\sitem[-]
  its queries at position~$w=x_1 \cdots x_n$ are pairs~$(i,q)$
  consisting of an index~$1\leq i\leq n$ and a query~$q$ at position~$x_i$
  in the game~$A$.
\sitemize
The polarities of moves and queries are inherited from the game~$A$
in the expected way, and the residual relation is defined as
for the tensor product.
%
%
Interestingly, the resulting multi-bracketed game~$! A$
defines the free commutative comonoid associated
to the well-bracketed game~$A$ in the category~$\ses$.
Hence, the category~$\ses$ defines a 
model of multiplicative exponential linear logic.
This model is \emph{degenerate} in the sense
that the tensor product is equal to its dual, 
i.e. $\dual{(A\tensor B)} = \dual{A}\tensor\dual{B}$.
%

%

\paragraph{Fixpoints.}

The exponential modality together with the traced symmetric monoidal
structure on $\ses$ defines a fixpoint operator in $\ses$ as shown by
Hasegawa Masahito in~\cite{hasegawa:phd}. Remark that this
construction does not require that the category~$\ses$ is cartesian.

\section{A game model with resources}\label{sct:game}

We would like to construct a model of tensorial logic based on
negative multi-bracketed games.  However it is meaningless to
construct an affine modality on the category~$\negative$ itself
because its unit~$1$ is already a terminal object in the category.
So we need to introduce the notion of \emph{pointed game.}
\paragraph{Pointed games.}
A pointed game may be seen in two different ways:
(1) as a positive multi-bracketed Conway game,
with a unique initial Player move,
(2) as a negative multi-bracketed Conway game,
except that the hypothesis that there are no queries
at the root~$\ast$ is now relaxed for Player queries.
From now on, we adopt the first point of view, and thus see
a pointed game as a positive game with a unique initial move.
Now, a morphism~$\sigma: A\morph{} B$ in the category~$\ses$
is called \emph{transverse} when, for every play $mn$ of length 2
in the strategy~$\sigma:\dual{A}\tensor B$, the Opponent move~$m$
is in~$A$ and the Player move~$n$ is in~$B$.
We note $\pointed$ the subcategory of~$\ses$
with pointed games as objects, and well-bracketed transverse strategies
as morphisms.

\paragraph{{\bf Coalesced tensor.}}

Given $A,B \in \pointed$, the coalesced tensor $A \synchTensor B$
is the pointed game obtained from~$A\tensor B$ by synchronising
the two initial Player moves of $A$ and~$B$.
Remark that the coalesced tensor product preserves affine games,
and coincides there with the tensor product of~$\negative$.
The category $\pointed$ equipped with~$\synchTensor$
is symmetric monoidal.
It is not monoidal closed, but admits a tensorial negation. 
Besides, it inherits a trace operator from the category~$\ses$, which
is partial, but sufficient to interpret a linear PCF with resource
modalities. 

\paragraph{Tensorial negation.}
The negation~$\lnot A$ of a pointed game~$A$ is the pointed game
obtained by lifting the dual game~$\dual A$ with a Proponent move~$m$
which initiates one query. Then, every initial Opponent move in~$\dual A$
complies with this query.

\paragraph{Affine modality.} 
A pointed game~$A$ is called \emph{affine} when its unique
initial Player move does not initiate any query.
Note that $\negative$ is isomorphic to the full subcategory of affine
games in the category~$\pointed$.
The affine game $\affine A$ associated to a pointed game~$A$
is defined by removing all the queries initiated by the first move
--- as well as their residuals.
This defines an affine resource modality on~$\pointed$.
%
%
%

\paragraph{Exponential modality.} The exponential modality
$\exponential{}$ on pointed games is obtained by composing the two
adjunctions underlying the comonads~$\affine$ and~$!$
(defined in~$\S \ref{sct:conway games}$).
$$
\xymatrix {\Mcat \ar@/^1pc/[rr]^{\subseteq} \ar@{}[rr]|-{\bot} &&
  \ar@/^1pc/[ll]^{!} \negative \ar@/^1pc/[rr]^{\subseteq}
  \ar@{}[rr]|-{\bot}&& \pointed \ar@/^1pc/[ll]^{\affinespecial{}}}
$$
In particular, given a pointed game~$A$, $\exponential{A}$ is
defined as 
$$
\exponential{A} \quad \define \quad !\left(\affine{A}\right)
$$
%


%





\paragraph{Free coproducts.}

The category~$\pointed$ lacks coproducts to be a model of (full) tensorial logic.
We adjust this by constructing its free completion, noted $\fam \pointed$,
under small coproducts~\cite{abramsky-mccusker:families}.
Given a category $\C$, the objects of $\fam \C$ are families $\{A_i |i
\in I \}$ of objects of the category. A morphism from $\{A_i |i \in I \}$ to $\{B_j|j \in J
\}$ consists of a reindexing function $f: I \rightarrow J$ together
with a family of morphisms $\{ f_i: A_i \rightarrow B_{f(i)} | i\in I \}$
of the category~$\C$.

$\family$ is a pseudo-commutative monad on $\cat$ \cite{hyland-power}.
Hence, the $2$-monad for symmetric monoidal categories distributes
over $\family$.
Consequently, (1) the category~$\fam{\C}$ inherits the symmetric monoidal structure
of a symmetric monoidal category~$\C$, and (2) the coproduct of~$\fam{\C}$
distributes over that tensor product, and (3) $\family$ preserves monoidal adjunctions.
Besides, $\family$ preserves categories with finite products and
categories with a terminal object.
The construction thus preserves affine and exponential
modalities in the sense of $\S \ref{sct:catModel}$.
%
%
Gathering all those remarks, we obtain that:
\begin{prop}
$\fam{\pointed}$ is a model of tensorial logic.
\end{prop}
Moreover, the category~$\fam{\pointed}$ has a fixpoint operator
restricted on its singleton objects --- that is, objects~$\{A_i |i \in I \}$
where~$I$ is singleton.
This is sufficient to interpret a linear variant of the language PCF
equipped with affine and exponential resource modalities,
in the category $\fam{\pointed}$.
%
%



\section{Conclusion}\label{sct:conclusion}
In this paper, we integrate resource modalities in game semantics,
in just the same way as they are integrated in linear logic.
%
%
%
%
This is achieved by reunderstanding the very topography of the field.
More specifically, linear logic is relaxed into tensorial logic,
where the involutive negation
of linear logic is replaced by a tensorial negation.
Once this performed, it is possible to keep the best of linear logic:
resource modalities, etc. but transported in the language
of games and continuations.
Then, linear logic coincides with tensorial logic with
the additional axiom that the continuation monad
is commutative.
In that sense, tensorial logic is more primitive than linear logic,
in the same way that groups
are more primitive than abelian groups.
This opens a new horizon to the subject.
The whole point indeed is to understand in the future
how the theory of linear logic extends to this relaxed framework.
We illustrate this approach here by extending well-bracketing
in arena games to the full flavour of resources in linear logic,
using multi-bracketing.
%

%
%

%

\bibliographystyle{latex8}
\bibliography{resource_modalities.bib}

\section*{A sequent calculus for tensorial logic}
In the bilateral formulation of tensorial logic, the sequents are of
two forms: $\Gamma\vdash A$ where $\Gamma$ is a context, and~$A$ is a
formula; $\Gamma\vdash$ where $\Gamma$ is a context (the notation
$[A]$ expresses the unessential presence of $A$ in the sequent).
%
{\scriptsize
  \begin{center}
    \begin{tabular}{cc}
      \AxiomC{$\Gamma\vdash A$ \quad $\Delta\vdash B$}
      \RightLabel{$\tensor$-Right}
      \UnaryInfC{$\Gamma, \Delta \vdash A\tensor B$}
      \DisplayProof
      &
      \AxiomC{$\Gamma_1, A, B, \Gamma_2 \vdash [C] $}
      \RightLabel{$\tensor$-Left}
      \UnaryInfC{$\Gamma_1, A\tensor B, \Gamma_2 \vdash [C] $}
      \DisplayProof
      \\
      \\
      \AxiomC{}
      \RightLabel{Unit-Right}
      \UnaryInfC{$\vdash 1$}
      \DisplayProof
      &
      \AxiomC{$\Gamma\vdash [A]$}
      \RightLabel{Unit-Left}
      \UnaryInfC{$\Gamma, 1\vdash [A]$}
      \DisplayProof
      \\
      \\
      \AxiomC{$\Gamma, A \vdash $}
      \RightLabel{$\lnot$-Right}
      \UnaryInfC{$\Gamma \vdash \ \lnot A$}
      \DisplayProof
      &
      \AxiomC{$\Gamma\vdash A$}
      \RightLabel{$\lnot$-Left}
      \UnaryInfC{$\Gamma, \lnot A \vdash $}
      \DisplayProof
      \\
      \\
      \AxiomC{}
      \RightLabel{Axiom}
      \UnaryInfC{$ A \vdash A$}
      \DisplayProof
      &
      \AxiomC{$\Gamma \vdash A$ \quad $A,\Delta\vdash [B]$}
      \RightLabel{Cut}
      \UnaryInfC{$\Gamma , \Delta \vdash [B]$}
      \DisplayProof
      \\
      \\
      \begin{tabular}{c}
        \AxiomC{$\Gamma \vdash A$}
        \RightLabel{$\oplus$-Right-1}
        \UnaryInfC{$\Gamma \vdash A\oplus B$}
        \DisplayProof
        \\
        \\
        \AxiomC{$\Gamma \vdash B$}
        \RightLabel{$\oplus$-Right-2}
        \UnaryInfC{$\Gamma \vdash A\oplus B$}
        \DisplayProof
      \end{tabular}
      &
      \AxiomC{$\Gamma, A\vdash C$ \quad $\Gamma, B\vdash C$}
      \RightLabel{$\oplus$-Left}
      \UnaryInfC{$\Gamma, A\oplus B \vdash C$}
      \DisplayProof
      \\
      \\
      No right introduction rule for $0$
      &
      \AxiomC{}
      \RightLabel{$0$-Left}
      \UnaryInfC{$\Gamma, 0 \vdash A$}
      \DisplayProof
      \\
      \\
     \AxiomC{$\exponentialpetit \Gamma \vdash A$}
      \RightLabel{Strengthening}
      \UnaryInfC{$\exponentialpetit \Gamma \vdash \ \exponentialpetit A$}
      \DisplayProof
      &
      \AxiomC{$\Gamma, A \vdash [B]$}
      \RightLabel{Dereliction}
      \UnaryInfC{$\Gamma, \exponentialpetit A\vdash [B]$}
      \DisplayProof
      \\
      \\
      \AxiomC{$\Gamma \vdash [B]$}
      \RightLabel{Weakening}
      \UnaryInfC{$\Gamma, \exponentialpetit A \vdash [B]$}
      \DisplayProof
      &
      \AxiomC{$\Gamma, \exponentialpetit A, \exponentialpetit A \vdash [B]$}
      \RightLabel{Contraction}
      \UnaryInfC{$\Gamma, \exponentialpetit A\vdash [B]$}
      \DisplayProof
      \\
      \\
      \AxiomC{$\affinepetit\Gamma \vdash A$}
      \RightLabel{Strengthening}
      \UnaryInfC{$\affinepetit\Gamma \vdash \ \affinepetit A$}
      \DisplayProof
      &
      \AxiomC{$\Gamma, A \vdash [B]$}
      \RightLabel{Dereliction}
      \UnaryInfC{$\Gamma, \affinepetit A\vdash [B]$}
      \DisplayProof
      \\
      \\
      \AxiomC{$\Gamma \vdash [B]$}
      \RightLabel{Weakening}
      \UnaryInfC{$\Gamma, \affinepetit A \vdash [B]$}
      \DisplayProof
      &
    \end{tabular}
  \end{center}
}
\noindent
The monolateral formulation requires to polarize formulas,
and to clone each construct into a negative counterpart.
$$
\begin{array}{lp{.5em}r}
\mbox{Positives} && 0 \hspace{0.25ex} \ | \ \hspace{0.3ex}
1 \hspace{0.3ex} \ | \hspace{0.2ex} \downarrow L \ | \ \hspace{0.2ex}
P \tensor Q \hspace{0.2ex} \ | \ P \oplus Q \ | \ \affine P \ |
\ \exponential P \\ \mbox{Negatives} && \bot \ | \ \top \ | \uparrow P
\ | \ L \Par M \ | \ \hspace{0.25ex} L \with M \hspace{0.25ex} \ |
\ \affinequestionmark L \ | \ \questionmark L
\end{array}
$$
It is then possible to reformulate all the sequent above,
as illustrated below by the right and left introduction of $\tensor$.
{\scriptsize
  \begin{center}
    \begin{tabular}{cc}
      \AxiomC{$\vdash \Gamma, P$ \quad $\vdash \Delta, Q$}
      \RightLabel{$\tensor$ \mbox{\textcolor{white}{-Right}}}
      \UnaryInfC{$\vdash \Gamma, \Delta, P\tensor Q$}
      \DisplayProof
      &
      \AxiomC{$\vdash \Gamma_1, L, M , \Gamma_2 , [P]$}
      \RightLabel{$\Par$ \mbox{\textcolor{white}{-Left}}}
      \UnaryInfC{$\vdash \Gamma_1, L\Par M, \Gamma_2 , [P] $}
      \DisplayProof
    \end{tabular}
  \end{center}
}

\end{document}